\documentclass[12pt,reqno,fleqn]{amsart}
\usepackage{amsmath,amstext,amsthm,amssymb,amsxtra}
\usepackage[colorlinks,citecolor=red,pagebackref,hypertexnames=false]{hyperref}
\usepackage{pgf,tikz}
\usepackage{hyperref}
\usepackage{pdfsync}
\allowdisplaybreaks

\newcommand{\RR}{{\mathbb{R}}}
\newcommand{\NN}{{\mathbb{N}}}

\newcommand{\eps}{\varepsilon}
\newcommand{\bp}{\noindent {\it Proof}.\,\,}
\newcommand{\ep}{\hfill$\Box$ \vskip 0.08in}

\newtheorem{proposition}{Proposition}[section]
\newtheorem{theorem}[proposition]{Theorem}
\newtheorem{lemma}[proposition]{Lemma}

\newtheorem{definition}{Definition}[section]

\usepackage{color}
\usepackage{enumerate}
\usepackage{enumitem}

\definecolor{gr}{rgb}   {0.,   0.8,   0. }
\definecolor{bl}{rgb}   {0.,   0.1,   1. }
\definecolor{mg}{rgb}   {0.7,  0.,    0.7}

\begin{document}

\author{Koushik Ramachandran}

\thanks{The author was partially supported by the NSF grant DMS - 1067886}


\address{Department of Mathematics, Purdue University, 150 N. University Street, West Lafayette, IN 47907-2067, USA}
\email{kramacha@math.purdue.edu}

\title[Harmonic functions in certain domains]{Asymptotic behavior of positive harmonic functions in certain unbounded domains}
\date{\today}


\maketitle

\vspace{0.1in}

\begin{abstract}
We derive asymptotic estimates at infinity for positive harmonic
functions in a large class of non-smooth unbounded domains. These
include domains whose sections, after rescaling, resemble a Lipschitz
cylinder or a Lipschitz cone, e.g., various paraboloids and horns.
\end{abstract}

\section{Introduction}
\noindent This paper involves a study of the growth of positive harmonic functions. Let $\Omega$ be an unbounded Lipschitz domain in $\mathbb{R}^d.$ We use the notation $(x_1, Y)$ to denote a point in $\mathbb{R}^d.$ Here $x_1\in\mathbb{R}$ and $Y\in\mathbb{R}^{d-1}.$ Sometimes we will also use spherical coordinates $\xi= (r,\omega),$ where $r=|\xi|$ and $\omega =\frac{\xi}{r}.$ Let $u$ be a positive harmonic function in $\Omega$ and suppose that $u=0$ on the boundary $\partial\Omega.$ We are interested in the asymptotic behavior of $u$ at infinity. Depending on the nature of the domain, we consider either of the quantities
\begin{equation}
M(r)=\sup_{\{Y:(r,Y)\in\Omega\}}u(r,Y)
\end{equation}

\begin{equation}
\widetilde M(r)=\sup_{\{\omega:(r,\omega)\in\Omega\}}u(r,\omega),
\end{equation}
\noindent and obtain asymptotics of these quantities as $r$ tends to infinity. The problem of studying the growth of positive harmonic functions is a very old one. Friedland and Hayman, in their paper, (\cite{Fr-Hay}, Corollary to Theorem D) attribute a lower bound on $\widetilde M$ to Huber. Their lower bound is the following.
\begin{equation}\label{eqHF}
\widetilde M(r)\geq C\exp\left(\int_{e}^{r/2}\frac{\alpha_1(t)}{t}dt\right),
\end{equation}
\noindent where $\alpha_1(t)$ is the so called characteristic constant of $\Omega\cap\{\xi: |\xi|=t\}.$ For the definition of the characteristic constant see section $3.1$. The above lower bound is valid in all domains and in general, one cannot infer anything more about the growth. Huber's result as well as results of Hayman and Friedland, both are based on a differential inequality technique attributed to Carleman in dimension $d=2.$

\vspace{0.2in}

\noindent The main purpose of our paper is to show that the lower estimate \eqref{eqHF} becomes an asymptotic equality for certain classes of domains that we introduce in the following paragraphs.

\vspace{0.2in}

\noindent In this paper we deal only with dimension $d\geq3.$ When $d=2,$ one can use conformal mapping techniques to find an asymptotic equality for $M(r).$ This method is due to Warschawski and is known in the literature as the Ahlfors-Warschawski estimates. See reference (\cite{War}). For $d\geq 3,$ conformal mapping methods are not applicable.

\vspace{0.2in}

\noindent A problem that is closely related to that of obtaining asymptotics is the estimation of harmonic measure near infinity. Indeed, one can obtain a lower bound for harmonic measure if the growth of $\widetilde M(r)$  is known. Let us show how this is done. Assume that $\widetilde M(r)$ is known. Let $S_r = \Omega\cap\{\xi: |\xi|=r\}.$ Fix $x_0\in\Omega.$ For all $r$ large enough, we will have $|x_0|<r.$ For such $r,$ denote by $\omega_{x_0}(r)$ the harmonic measure of $S_r$ with respect to the domain $ \Omega\cap\{\xi: |\xi|<r\},$ based at the point $x_0.$ We can assume without loss of generality that $u(x_0)=1.$ Then, using the definition of harmonic measure along with the fact that $u=0$ on $\partial\Omega,$ we have
\begin{equation}\label{Hmeasure}
1=u(x_0)=\int_{S_r}u(\xi)d\omega_{x_0}(\xi)\leq \widetilde M(r)\omega_{x_0}(r).
\end{equation}
From equation \eqref{Hmeasure}, we have
\begin{equation}\label{eqrecip}
\omega_{x_0}(r)\geq\frac{1}{\widetilde M(r)}.
\end{equation}

\noindent If $\widetilde M(r)$ is known, \eqref{eqrecip} gives a lower bound on harmonic measure. On the other hand, obtaining an upper bound for harmonic measure is not simple, but in domains with a nice geometrical shape the Carleman method c.f. \cite{Hal}, will give us good upper bounds.  In cones, harmonic measure estimates near infinity have been obtained by Burkholder \cite{Burk}, Essen and Haliste \cite{Ess-Hal} and some others who we list in the references. Banuelos and Carroll \cite{Ban-Carr}, established harmonic measure estimates in paraboloid-type domains
$$P_{\alpha}= \{(x,Y)\in \mathbb{R}\times\mathbb{R}^{d-1}: \,x_1 > 0, \,|Y|< Ax_1^{\alpha}\},$$
where $A>0$ and $0<\alpha<1.$ In proving their result, they used the estimates of Warschawski \cite{War} in dimension $d=2,$ and for $d\geq 3$ they used the results in Hayman and Carroll \cite{Car-Hay}. The literature that deals with estimating harmonic measure is vast. Quite a few of these results have been proved using methods of probability. This is not surprising considering the relation between Brownian motion and harmonic measure.
\vspace{0.2in}

\noindent Going back to our problem, the asymptotics of $u$ are well known when $\Omega$ is a cylinder or a cone, for in these cases an explicit form of a solution $u$ is known, cf. Ancona \cite{Ancona}. It is then easy to read off the exact growth of $M(r),$ or correspondingly $\widetilde M(r),$ from these formulas. Usually, in most other domains, explicit formulas for a solution $u$ are unavailable. Despite this, one would think that finding sharp asymptotics in reasonably simple domains should not be difficult. Somewhat surprisingly, sharp asymptotics are unknown even in the case of a paraboloid.

\vspace{0.1in}

\noindent Cranston and Li \cite{Cran-Li}, obtained pointwise upper and lower bounds for positive harmonic functions vanishing on the lateral side of horn-shaped domains of the form
$$H_f = \{(x_1, Y)\in \mathbb{R}\times\mathbb{R}^{d-1}:\,x_1>0,\,|Y|<f(x_1)\},$$
where $f:(0,\infty)\rightarrow (0,\infty)$ and in addition satisfies some other technical conditions which seem difficult to verify. Their proof is also considerably more complicated than what we present here. In his paper \cite{DeBlassie}, DeBlassie derived asymptotics for $u$ at infinity, in domains given by $$\Omega_a=\{(x_1, Y)\in \mathbb{R}\times\mathbb{R}^{d-1}:\,x>0,\,|Y|<a(x_1)\},$$
where $a$ is a positive Lipschitz function defined on $[0,\infty)$, satisfying certain growth constraints. The difference in the domains $H_f$ and $\Omega_a$ arise because of the different assumptions made on $f$ and $a.$  The class of domains $\Omega_a$ include both the paraboloid-type domains $P_{\alpha}$ and certain type of horn-shaped domains $H_f.$ The domains that were studied by DeBlassie, and Banuelos and Carroll, are all rotationally invariant in the following sense : if $(x_1, Y)\in\Omega,$ then it is also true that $(x_1, \tilde Y)\in\Omega,$ whenever $|Y|=|\tilde Y|.$ This fact is very much used by the authors, in proving their results, and the methods do not extend to domains that are not rotationally invariant. Therefore the question arises as to how to obtain asymptotics in such domains.

\vspace{0.1in}

\noindent In this paper, we introduce two classes of domains, namely, cylinder-like domains and cone-like domains and establish asymptotics in them. These classes of domains include many domains which are not rotationally invariant. It is also reassuring that they contain a number of simple domains that are easy to describe, including those studied by DeBlassie. Indeed, the simplest examples of cylinder-like domains are the paraboloid-type domains $P_{\alpha}.$ Roughly speaking, a cylinder-like domain is one which tends to a cylinder after a translation, and scaling by an appropriate Lipschitz function. Along similar lines, a cone-like domain is one which tends to a cone after scaling by a well-chosen factor.

\vspace{0.1in}

\noindent We will now briefly explain our method of obtaining asymptotics. For simplicity we give the method in the case of cone-like domains. Suppose that we want to find asymptotics of a positive harmonic function $u,$ in a cone-like domain $\Omega.$ Our first step is to scale portions of $\Omega$ that are far away. By definition, these scaled domains $\Gamma_n$ will tend to a cone $C.$ Using the scaling map, we define a sequence of positive harmonic functions $\{v_n\}$ which live on $\Gamma_n$. The $v_n$ strongly reflect the behavior of $u$ near infinity. The next step is to prove that one can extract a subsequence of $\{v_n\}$ which converges uniformly on compact subsets of $C$, to a harmonic function $v$. Since $\Gamma_n$ tend to a cone, it is natural to guess that the support of $v$ lies in the cone. We prove that this is indeed the case. This part of the proof makes use of a Boundary Harnack type estimate and also Holder estimates near the boundary. Let us continue with the proof. As mentioned before behavior of harmonic functions in cones is completely known. Finally, using elementary analysis we infer the asymptotics of $u$ from the known behavior of $v.$ We employ a similar procedure to find asymptotics in cylinder-like domains.

\vspace{0.1in}

\setcounter{equation}{0}
\section{Cylinder-like Domains}
\subsection{Notation and definitions}
\vspace{0.1in}

\noindent In this section we denote points in $\mathbb{R}^d$ by $(x_1, Y),$ where $x_1\in\mathbb{R}$ and $Y\in\mathbb{R}^{d-1}.$ We say that a function $\Phi : \mathbb{R}^{d-1}\rightarrow\mathbb{R}$ is a Lipschitz function if there exists $k < \infty$ such that
$$|\Phi(Y) - \Phi(\hat Y)|\leq k|Y - \hat Y|$$
for all $Y,$ $\hat Y$ in $\mathbb{R}^{d-1}.$

\begin{definition}\label{Lipschitz}
A domain $\Omega$ (bounded or unbounded) is called an uniform Lipschitz domain if there exists a $k < \infty,$ and for every $p\in\partial\Omega$ there is
\begin{enumerate}
\item A neighborhood $U_p$ of $p$
\item a Lipschitz function $\Phi_p$ with Lipschitz constant not greater than $k$,
\end{enumerate}
such that $\Omega\cap U_p = \{(x_1, Y): y_{d-1} > \Phi_p (x_1, y_1, y_2...y_{d-2})\}\cap U_p,$ in some orthonormal coordinate system which may depend on the point $p.$
\end{definition}

We will denote a unbounded uniform Lipschitz domain by UUL domain. We write $e_1$ for the point in $\mathbb{R}^d$ with coordinates $(1,0,...0).$ For $t>0$, we set
\begin{equation}\label{eqbig}
\Omega_t=\Omega\cap\{\ t/2 < x < 3t/2 \}
\end{equation}

\begin{equation}\label{eqgam}
\Gamma_t = (\Omega_t-te_1)/a(t).
\end{equation}

\noindent A Lipschitz cylinder will mean a set of the form $\mathbb{R}\times D$, where $D$ is a bounded Lipschitz domain in $\mathbb{R}^{d-1}.$ By $d_H(A_1,A_2)$ we mean the Hausdorff distance between the sets $A_1$ and $A_2.$

\vspace{0.1in}

\begin{definition}\label{defcylinder}
An UUL domain $\Omega$ is called cylinder-like, if there exists a positive Lipschitz function $a$ defined on the interval $[0,\infty)$ and a family of Lipschitz cylinders $\mathfrak C_t=\mathbb{R}\times D_t$ such that for every compact $K$\\
$$ d_H(\Gamma_t\cap K, \mathfrak C_t\cap K)\to 0\hspace{0.1in}\mbox{as}\hspace{0.1in} t\to\infty.$$

\noindent We will also impose the following conditions:

\begin{enumerate}
\item $ a'(t)\searrow 0$ as $t\to\infty.$
\item  The Lipschitz constants of the cylinders $\mathfrak C_t$ are uniformly bounded above and below.
\item All the cylinders $\mathfrak C_t$ contain a fixed Lipschitz cylinder $L_0$ and are themselves contained in another Lipschitz cylinder $L_1$.
\item If $\lambda(t)$ denotes the principal eigenvalue of the $\Delta_{d-1}$ acting on the domain $D_t,$ then $\lambda(t)$ is continuous as a function of $t$.
\end{enumerate}
\end{definition}

\vspace{0.1in}

\noindent We make a few remarks at this point. Condition $(3)$ in definition \ref{defcylinder} ensures that there exist positive constants $\lambda_0,\lambda_1$ such that $\lambda_0\leq\lambda(t)\leq\lambda_1$. Another fact which we will use later and which we want to point out is that the condition $a'(t)\searrow 0$ implies, by the mean value theorem, that
\begin{equation}\label{abyt}
\frac{a(t)}{t}\to 0 \hspace{0.1in}\mbox{as}\hspace{0.1in} t\to\infty.
\end{equation}

\subsection{Theorem for cylinder-like domains}
\setcounter{equation}{3}

\vspace{0.1in}

\begin{theorem} \label{cylinders}

Let $u$ be a positive harmonic function on a cylinder-like domain $\Omega$. Suppose that $u=0$ on $\partial\Omega$. For $t>0,$ denote $M(t)=\sup_{\{(x,Y)\in\Omega: x=t\}}u(x,Y).$ Then, $M(t)$ satisfies the following asymptotic formula.
$$ \log M(t)=(1+\rho(t))\int_{1}^{t}\dfrac{\sqrt{\lambda(\tau)}}{a(\tau)}d\tau,\hspace{0.1in}\mbox{where}\hspace{0.1in} \rho(t)\to 0\hspace{0.05in}\mbox{as}\hspace{0.05in} t\to\infty.$$

\end{theorem}

\bp Fix a positive sequence $\{t_n\}$, $t_n\to\infty$. Let $\Omega_{t_n}$ and $\Gamma_{t_n}$ be as defined in equations \eqref{eqbig} and \eqref{eqgam} respectively. By hypothesis we have that on every compact $K$, $d_H(\Gamma_{t_n}\cap K, \mathfrak C_{t_n}\cap K)\to 0\hspace{0.1in}\mbox{as}\hspace{0.1in} n\to\infty.$ Now we use the fact that the set of Lipschitz cylinders is compact in the topology of Hausdorff distance, cf. \cite{Buc-But}. Hence there exists a subsequence $\{t_{n_k}\}$ such that $\{\mathfrak C_{t_{n_k}}\}$ converges to an open set $\mathfrak C$ in the Hausdorff distance. It is evident that the limit $\mathfrak C$ is also a Lipschitz cylinder. Hence $\mathfrak C$ has the form $\mathbb{R}\times D$ for some bounded Lipschitz domain $D\subset\mathbb{R}^{d-1}.$ Applying the triangle inequality we can now conclude that there exists a subsequence $\{t_{n_k}\}$ such that

\begin{equation}\label{eqHlimit}
d_H(\Gamma_{t_{n_k}}\cap K, \mathfrak C\cap K)\to 0\hspace{0.1in}\mbox{as $k\to\infty$, for every compact K.}
\end{equation}

\noindent To keep the notations simple we will rename the subsequence $\{t_{n_k}\}$ to be $\{t_n\}$, and the domains $\Gamma_{t_{n_k}}$ and ${\Omega}_{t_{n_k}}$ to be $\Gamma_n$ and $\Omega_n$ respectively. Having fixed such a subsequence we now observe that the mapping $z\mapsto a(t_n)z+t_ne_1$ maps $\overline{\Gamma}_n$ into $\overline{\Omega}_n.$ Furthermore, for $z\in\Gamma_n$, define $$ v_n(z)=\frac{u(a(t_n)z+t_ne_1)}{M(t_n)}.$$ Note that $v_n$ satisfies $v_n(0)\leq1$ and

\begin{equation}\label{harm}\Delta v_n=\dfrac{(a(t_n))^2}{M(t_n)}\Delta u=0,\hspace{0.1in}\mbox{in $\Gamma_n$.}
\end{equation}

\noindent Hence each $v_n$ is a positive harmonic function on $\Gamma_n$ with values $0$ on the lateral part of $\partial\Gamma_n.$ For the following we need a lemma.

\begin{lemma}\label{lemmageneral}
Let $\Lambda\in\RR^d$ be a domain, and suppose that $\{\Lambda_n\}_n$ is a sequence of domains such that for every compact $K,$
$$d_H(\Lambda_n\cap K, \Lambda\cap K)\to 0 \quad\mbox{as}\quad n\to\infty,$$
and all the $\Lambda_n$ contain a common point $\xi_0.$ Let $\{h_n\}$ be a sequence of positive harmonic functions on $\Lambda_n$, $n\in\NN$. Suppose further that the sequence $\{h_n\}_n$ is bounded at the point $\xi_0$, that is, there exists a constant $M>0$ such that
\begin{equation} \label{eq1} h_n(\xi_0)<M\quad \mbox{for all}\quad n\in\NN.\end{equation}
Then, there exists a subsequence $\{h_{n_{i}}\}\subset\{h_n\}_n$ which converges uniformly on compact subsets of $\Lambda$.
\begin{proof} This lemma is, in fact, a mild modification of the classical convergence theorem for harmonic functions in a domain following from interior estimates (cf., e.g., \cite{GT}, Theorem~2.11). Indeed, consider some sequence of $\{K_j\}_j$, compacta in $\Lambda$, such that
 $$\Lambda=\bigcup_{j\geq 1}\ K_j, \quad K_j\subset K_{j+1}^o,\quad j\in\NN,$$
where $K_j^o$ denotes the interior of $K_j.$ Without loss of generality we can assume that $\xi_0\in K_j$ for all $j\in\NN$, extending $K_j$ to contain $\xi_0$, if necessary. Since $\Lambda_n\to\Lambda$ in the Hausdorff distance on every compact, we have for every $j$ there exists $N(j)$ such that for all $n\geq N(j)$ we have $K_j\subset\Lambda_n$, and hence, $h_n$ are harmonic on $K_j$ for all such $n$.

Due to the interior Harnack inequality and \eqref{eq1}, for every $j$ there exists $C_j$ such that
\begin{equation}\label{eq2}\sup_{z\in K_j} h_n(z)\leq C_j \inf_{z\in K_j} h_n(z)\leq C_j M,\quad \mbox{for all } n\geq N(j),\end{equation}
and moreover, for any multi-index $\alpha$ there exists $C_{j,\alpha}$ such that
\begin{equation}\label{eq3}
\sup_{z\in K_j} |\nabla^{\alpha}h_n(z)|\leq C_{j,\alpha}\sup_{z\in K_{j+1}} h_n(z) \leq C(j,\alpha,M), \quad \mbox{for all } n\geq N(j+1).
\end{equation}
The first inequality in \eqref{eq3} follows from the standard interior estimates (cf., e.g., \cite{GT}, \S{2.7}), and the second one is a consequence of \eqref{eq2}. Therefore the sequence $\{h_n\}_{n\geq N(j+1)}$ is equicontinuous and uniformly bounded on $K_j$, and hence, by Arcela-Ascoli Theorem, we can extract a subsequence $\{h_n^j\}$ which converges uniformly on $K_j$. Clearly, one can choose $\{h_n^j\}$ to be a subsequence of $h_n^{j-1}$ for $j\geq 2$.
We may now employ Cantor's diagonal procedure to pick a subsequence $\{h_{n_i}\}$ which converges uniformly on each of the $K_j$.  Finally, since every compact $K\subset\Lambda$ is a subset of some $K_j$, the sequence $\{h_{n_i}\}$ converges uniformly on $K$, as desired.
\end{proof}
\end{lemma}

Observe that $v_{t_n}(0)\leq 1$, so we see that the conditions of the above lemma are satisfied with $v_n$ in place of $h_n$ and $\Gamma_n$ in place of $\Lambda_n.$ Applying the lemma we deduce that a subsequence, which we still denote $ v_n$, converges to a function $v$ uniformly on compact subsets of the cylinder ${\mathfrak C}.$ Since $v$ is the uniform limit of harmonic functions, it is harmonic. We now claim that $v=0$ on $\partial {\mathfrak C}.$

\begin{lemma}\label{lemma5} Let $v$ be a harmonic function on the Lipschitz cylinder ${\mathfrak C}$ given by
\begin{equation}\label{eqDefv}
v=\lim_{n\to\infty} v_n, \quad v_n(z)=\dfrac{u(a(t_n)z+t_ne_1)}{M(t_n)},\quad \mbox{for $z\in \Gamma_n$},
\end{equation}
where the convergence is uniform on compact subsets of $\mathfrak C$ and $v_n$ vanishes on the lateral part of $\partial\Gamma_n.$ Then $v$ extends to a continuous function on $\overline{{\mathfrak C}}$ with $v|_{\partial {\mathfrak C}}=0$ and furthermore, for each $r\in\mathbb{R},$
$$\sup_{\{x=r\}\cap\mathfrak C}{v_n}\to\sup_{\{x=r\}\cap\mathfrak C}{v}.$$
\end{lemma}

\vspace{0.1in}

\bp Fix real numbers $a < b$ and consider the set $G = \{a< x_1 < b\}\cap\mathfrak C.$ Let $\epsilon > 0$ be given. Now $\overline G$ is a compact set. Therefore by hypothesis there exists $N = N(\epsilon, G)$ such that for $n\geq N$,
\begin{equation}\label{bdrydis}
d_H(\Gamma_n\cap\overline G, \mathfrak C\cap\overline G)<\epsilon.
\end{equation}

\noindent Furthermore, by taking $N$ larger if necessary, we may assume that for $n\geq N$ only the lateral portion $\partial\Gamma_n^L$ of $\partial\Gamma_n$ intersects $\overline G.$ The Lipschitz cylinder $\mathfrak C$ has the following property : there exists a fixed cone $C(\theta, h),$ of opening $\theta$ and height $h$ such that at every boundary point of $\mathfrak C$ one can draw an exterior cone congruent to $C(\theta, h).$ Since $d_H(\Gamma_n\cap K, \mathfrak C\cap K)\to 0$ on every compact $K$, $\Gamma_n$ also has the same property with the same $C(\theta, h)$ provided $n$ is large. In particular there exists a fixed $r_0 > 0$ independent of $n$ such that for $n\geq N$ and $z\in\partial\Gamma_n^L$, $B_{r_0}(z)\cap\Gamma_n$ can be represented as the region above the graph of some Lipschitz function, with the corresponding Lipschitz constants being uniformly bounded from above and below by positive constants independent of $n.$

\vspace{0.1in}

\noindent Let us take any $\epsilon,$ $0 < \epsilon < r_0/100,$ with $r_0$ being as in the preceding paragraph. Let $N$ be the corresponding $N(\epsilon, G)$ for which \eqref{bdrydis} holds. Let $\partial G^L$ be the lateral portion of $\partial G.$ Denote $(\partial G^L)_{\epsilon}= \{x: d(x, \partial G^L) < \epsilon\}.$ Then
\begin{equation}\label{bdryclose}
n\geq N\Rightarrow \partial\Gamma_n^L\cap \{a < x_1 < b\}\subset (\partial G^L)_{\epsilon}.
\end{equation}

\noindent Next cover $(\partial G^L)_{\epsilon}$ by finitely many balls centered at $\partial G^L,$ i.e,
\begin{equation}\label{cover}
(\partial G^L)_{\epsilon}\subset \bigcup\limits_{i=1}^m{B(z_i, \delta)},
\end{equation}

\vspace{-0.1in}

\noindent where $z_i\in\partial G^L$ and $\delta = r_0/2.$ This is possible by our choice of $\epsilon$ and $\delta.$ Now take any $x\in G\cap B(z_1, \epsilon).$ For $n\geq N$, if $x\notin\Gamma_n$ extend the definition by setting $v_n(x) = 0.$ Otherwise $x\in\Gamma_n,$ and using $x\in B(z_1, \epsilon)$ along with \eqref{bdryclose}, we will get
\begin{equation}\label{dtolateral}
d(x, \partial\Gamma_n^L) < 2\epsilon.
\end{equation}

\noindent Furthermore $v_n$ is a positive harmonic function and $v_n = 0$ on the Lipschitz portion $ \partial\Gamma_n^L.$ Therefore, by H\"{o}lder estimates near the boundary $\partial\Gamma_n^L$ (c.f \cite{KenigBook}), there exist constants $C>0$ and $0 < \alpha\leq 1$ such that
\begin{equation}\label{Holder}
v_n(x)\leq C\left(\frac{d(x, \partial\Gamma_n^L)}{\delta/2}\right)^{\alpha}\sup_{B(z_1,\delta)\cap\Gamma_n}{v_n}.
\end{equation}

\noindent In general the constants $C$ and $\alpha$ in equation \eqref{Holder} will depend on the Lipschitz character of the domains, however by our remark earlier, $\partial\Gamma_n^L$ have Lipschitz constants uniformly bounded from above and below and hence we may choose $C$ and $\alpha$ independent of $n.$ Next, we can control the supremum in \eqref{Holder} by the Carleson estimate c.f \cite{Aikawa},
\begin{equation}\label{Carleson}
\sup_{B(z_1,\delta)\cap\Gamma_n}{v_n}\leq C_1v_n(z_{1,\delta}),
\end{equation}
where $ z_{1,\delta}$ is a point in $B(z_1, \delta)\cap\Gamma_n\cap G$ with $|z_1 - z_{1,\delta}|\geq 3\delta/4= 3r_0/8.$ Again, $C_1$ in equation \eqref{Carleson} will usually depend on $n$ but using uniform Lipschitz condition of $\Gamma_n$, we may choose it independent of $n.$ Finally let $G_{r_0/8}$ be the compact set $\{x\in G: d(x, \partial G)\geq r_0/8\}.$ Then using the fact that $\{v_n\}$ converges uniformly on compacts to $v,$ we will have for $n\geq N$ that $v_n$ are all uniformly bounded in $G_{r_0/8}$ by a constant $C_2$ that is independent of $n.$ In particular since $z_{1,\delta}\in G_{r_0/8},$ we have for $n\geq N$
\begin{equation}\label{unifbd}
v_n(z_{1,\delta})\leq C_2.
\end{equation}

\noindent Combining equations \eqref{dtolateral}, \eqref{Holder}, \eqref{Carleson} and \eqref{unifbd}, we get
\begin{equation}\label{small}
v_n(x)\leq C_3\epsilon^{\alpha},
\end{equation}

\noindent for a constant $C_3$ independent of $n.$ A similar estimate as \eqref{small} holds for points $\epsilon$ close to $z_i$ in the other balls $B(z_i, \delta)$ for $2\leq i\leq m$ with the same constant $C_3.$ Since $\bigcup\limits_{i=1}^m{B(z_i, \delta)}$ covers $\partial G^L$, we can conclude that $v = \lim v_n$ satisfies $v = 0$ on $\partial G^L .$ But $\partial G^L = \{a < x_1 < b\}\cap\partial\mathfrak C$ (by the definition of $G$). Since $a$ and $b$ were arbitrary to start with, we have $v =0$ on $\partial\mathfrak C.$

\vspace{0.1in}

\noindent Moreover it is also true that $v_n$ converges to $v$ uniformly on each section of the form $\{x_1 = r\}\cap\mathfrak C.$ This is because, for points $x$ on $\{x_1 = r\}\cap\mathfrak C$ which are close to the boundary $\{x_1 =r\}\cap\partial\mathfrak C,$ $v_n$ is small by \eqref{small} and since we now know $v= 0$ on $\partial\mathfrak C,$ we can apply the same H\"{o}lder estimates as before to show $v(x)$ also satisfies an estimate of the form \eqref{small}. But this means we have $v_n$ converging uniformly to $v$ near the boundary $\{x_1 =r\}\cap\partial\mathfrak C.$ On the other hand, for points well inside $\{x_1 = r\}\cap\mathfrak C$ ( meaning a compact), we already know that $v_n$ converges uniformly to $v.$ Therefore $v_n$ converges uniformly to $v$ on the whole section $\{x_1 = r\}\cap\mathfrak C.$ From uniform convergence we get the convergence of the sup norms. That is,
\begin{equation}\label{Sup}
\sup_{\{x=r\}\cap\mathfrak C}{v_n}\to\sup_{\{x=r\}\cap\mathfrak C}{v}.
\end{equation}

\noindent This finishes the proof of the lemma.
\ep

Combining equation \eqref{harm} with Lemma \ref{lemma5} we deduce that the function $v$ defined by \eqref{eqDefv} satisfies
\begin{equation}\label{eqPropv}
\Delta v=0 \mbox{ in } {\mathfrak C}, \quad v|_{\partial {\mathfrak C}}=0,\quad v>0  \mbox{ in } {\mathfrak C}.
\end{equation}

\noindent Therefore, $v$ can be represented as
\begin{equation}\label{eq16}
v(z)=A_1 e^{\sqrt {\tilde\lambda} z_1}\psi_1(z_2,...,z_d)+A_2 e^{-\sqrt{\tilde\lambda} z_1}\psi_1(z_2,...,z_d),\quad z\in {\mathfrak C},
\end{equation}

\noindent where $\tilde\lambda$ is the first eigenvalue of the Dirichlet Laplacian in the Lipschitz domain $D$ in $\RR^{d-1}$ and $\psi_1$ is the corresponding eigenfunction. For the above statement one can refer to \cite{Miyamoto}. Note that if $A_2\neq 0$, $v$ grows exponentially as $z_1$ approaches $-\infty$. If that is the case, then it follows that, in particular, for any $z_1$ sufficiently small  there exists $N=N(z_1)$ such that for all $n\geq N$ we have \begin{equation}\label{eqcontra}
u(a(t_n)z+t_ne_1) >M(t_n),
\end{equation}
for $z=(z_1,\tilde z)$ with $z_1<0$. However if $z_1<0$, then using the fact that the function $a$ is positive, we see that the $x$ coordinate of $a(t_n)z+t_ne_1 $ is less than $t_n.$ Recalling that $M(t)$ is the supremum of $u$ on $\{x=t\}$, we can  then conclude that $\dfrac{u(a(t_n)z+t_ne_1)}{M(t_n)}\leq 1$ with the help of the maximum principle. This is a contradiction to \eqref{eqcontra}. Hence $A_2=0$ and
\begin{equation}\label{eqhc}
v(z)=A_1 e^{\sqrt{\tilde\lambda} z_1}\psi_1(z_2,...,z_d),\quad z\in {\mathfrak C}.
\end{equation}
Let us now go back to the limit in equation \eqref{Sup}. Recalling the definition of $v_n$ we observe that the supremum on the left hand side of equation \eqref{Sup} is $\dfrac{ M(a(t_n)r+t_n)}{M(t_n)}.$ From equation \eqref{eqhc} above, we get that the supremum on the right hand side is $Ae^{\sqrt{\tilde\lambda} r}$ for a certain positive constant $A.$ Therefore, in the limit we get for every $r\in\mathbb{R}$,
$$\frac{ M(a(t_n)r+t_n)}{M(t_n)}\to Ae^{\sqrt{\tilde\lambda} r}\hspace{0.1in}\mbox{as}\hspace{0.1in} n\to\infty .$$

\noindent Plugging in the value 0 for $r$ we infer that $A=1$. We have now obtained the following.
$$\frac{ M(a(t_n)r+t_n)}{M(t_n)}\to e^{\sqrt{\tilde\lambda} r}\hspace{0.1in}\mbox{as}\hspace{0.1in} n\to\infty .$$

\noindent Taking logarithms and setting $f(t)=\log M(t)$, we get that for every sequence $\{t_n\}$, $t_n\to\infty$  there exists a subsequence $\{t_{n_k}\}$ such that for every $r\in\mathbb{R}$ we have,

\begin{equation}
\label{eq18} f(t_{n_k}+ ra(t_{n_k}))- f(t_{n_k})\to \sqrt{\tilde\lambda}r\hspace{0.1in}\mbox{as}\hspace{0.1in} n\to\infty.
\end{equation}

\noindent At this point we pause to remind that $\tilde\lambda= \lim_{k \to\infty}\lambda(t_{n_k}).$ Since equation \eqref{eq18} holds on  all subsequences, we may now conclude that,
\begin{equation}\label{eq18m}
f(t+ra(t))-f(t)-r\sqrt {\lambda(t)}\to 0\hspace{0.1in}\mbox{as}\hspace{0.1in} t\to\infty.
\end{equation}
Continuing with the proof of theorem \ref{cylinders}, set
\begin{equation}\label{phidef}
g_r(t)=t+ra(t), \quad t\in [0,\infty),\qquad \phi(t) =\int_{1}^{t}\frac{1}{a(y)}dy, \quad t\in [1,\infty).
\end{equation}

\noindent First observe that since the function $a$ is positive, $\phi$ is an increasing function. Therefore an inverse $\phi^{-1}$ exists. Also, by equation \eqref{abyt} in the introduction, $\frac{a(y)}{y}\to 0$ as $y\to\infty.$ Therefore if $y$ is large, $\frac{1}{y}\leq\frac{1}{a(y)}.$ Integration now yields
\begin{equation}\label{eqcomparison}
\log t-c_0\leq\int_{1}^{t}\frac{1}{a(y)}dy = \phi(t),
\end{equation}

\noindent for some constant $c_0.$ From equation \eqref{eqcomparison} we can conclude that $\phi(t)\to\infty$ as $t\to\infty.$ Since $\phi$ is increasing this also means that $\phi^{-1}(t)\to\infty$ as $t\to\infty.$ Hence if we set $s = \phi(t),$ we have that $s\to\infty\Rightarrow t\to\infty.$ Setting $h_r=\phi\circ g_r\circ\phi^{-1},$ we have
\begin {equation}
\label{eq19} f(t+ ra(t))= f\circ\ g_r(t)= f\circ\phi^{-1}(h_r\circ\phi(t)), \hspace{0.15in}\mbox{and } \hspace{0.15in} f(t)=f\circ\phi^{-1}(\phi(t)).
\end{equation}

\noindent Therefore using the notation
\begin{equation}\label{thetadef}
f_1:=f\circ \phi^{-1}, \hspace{0.04in}\theta:=\lambda\circ \phi^{-1},\hspace{0.02in}\mbox{and}\hspace{0.06in} s=\phi(t),
\end{equation}
along with equations \eqref{phidef} and \eqref{eq19}, we can now rewrite equation \eqref{eq18m} as
\begin{equation}
\label{eq20} f_1(h_r(s))-f_1(s)-r\sqrt {\theta(s)}\to 0,\hspace{0.1in}\ s\to\infty.
\end{equation}

\noindent Now,
\begin{eqnarray}\label{eq21}
h_r(s) &=& \phi\circ (g_r\circ\phi^{-1}(s))   \nonumber\\
                    &=& \int_{1}^{\phi^{-1}(s)+ra(\phi^{-1}(s))}\frac{1}{a(y)}dy  \nonumber\\
                    &=& \int_{1}^{\phi^{-1}(s)}\frac{1}{a(y)}dy+\int_{\phi^{-1}(s)}^{\phi^{-1}(s)+ra(\phi^{-1}(s))}\frac{1}{a(y)}dy \nonumber\\
                    &=& \phi\circ\phi^{-1}(s)+ \int_{\phi^{-1}(s)}^{\phi^{-1}(s)+ra(\phi^{-1}(s))}\frac{1}{a(y)}dy=s+r+\mu_r(s),
\end{eqnarray}

\noindent where $\mu_r(s)= \int_{\phi^{-1}(s)}^{\phi^{-1}(s)+ra(\phi^{-1}(s)}\frac{1}{a(y)}dy-r.$ With this in hand, \eqref{eq20} now becomes
\begin{equation}
\label{eqlim} f_1(s+r+\mu_r(s))-f_1(s)-r\sqrt{\theta(s)}\to 0,\hspace{0.1in}\ s\to\infty.
\end{equation}

\vspace{0.05in}

\noindent A crucial part of what follows is to estimate $\mu_r(s).$ But before we proceed towards this we make a remark which will help us later. We will be interested in the behavior of  $\frac{a(\xi_t)}{a(t)}$ as $t\to\infty$. Here $\xi_t$ denotes any point in the interval $[t,t+ra(t)]$ for a fixed $r$. We claim that the above ratio tends to $1$ as $t$ tends to infinity.
\begin{equation}\label{ratioclaim}
\frac{a(\xi_t)}{a(t)}\to 1 \hspace{0.02in}\mbox{as}\hspace{0.05in} t\to\infty.
\end{equation}

\noindent Indeed by the mean value theorem applied to the function $a$ we have,
\begin{equation}\label{eqratiolim}
\left|\frac{a(\xi_t)}{a(t)}-1\right| = \left|\frac{(\xi_t-t)a'(p_t)}{a(t)}\right|\leq |ra'(p_t)|,
\end{equation}
where $p_t$ is some point in the interval $[t,\xi_t].$ Since $a'\to\ 0$ we see that the claim is true. We are now ready to estimate $\mu_r(s)$.

\begin{lemma}\label{mu}
The function $\mu_r(s)$ defined above satisfies $\mu_r(s)\leq 0,$ and if $\epsilon>0$ is given then there exists a large number $S_{\epsilon}$ such that for $s > S_{\epsilon}$ and all $r\in\mathbb{R}$ we have
\begin{equation}
\label {eq23} |\mu_r(s)|\leq Cr^2\eps.
\end{equation}
\end{lemma}

\bp  Recall from equation \eqref{thetadef} that $\phi^{-1}(s)=t$ and that $t\to\infty$ as $s\to\infty$ by our remarks above. Now

\begin{eqnarray}
\label{eq24} \mu_r(s) &=& \int_{t}^{t+ra(t)}\frac{1}{a(y)}dy-r   \nonumber\\
                      &=& r\frac{a(t)}{a(t_r)}-r   \nonumber\\
                      &=& -r\frac{a(t_r)-a(t)}{a(t_r)}=-r\frac{(t_r-t).a'(\tilde t_r)}{a(t_r)}
\end{eqnarray}

\noindent In the above equations, $t_r$ denotes a point in the interval $[t, t+ra(t)]$ and $\tilde t_r$ a point in $[t, t_r]$. Note that we have used the mean value theorem for integrals in the second equality and the mean value theorem for functions in the fourth. Using the fact that $a$ and $a'$ are positive for large $t$ (the latter by condition ($1$) in definition \ref{defcylinder}), we can infer from \eqref{eq24} above that $\mu_r(s)\leq 0.$ We next notice that $|t_r-t|\leq |r|a(t)$ and from \eqref{ratioclaim} above $\frac{a(t)}{a(t_r)}$ tends to $1$ and hence is bounded above by a constant $C.$ Also because $a'\to 0$ near infinity, there is a number $S_{\epsilon}$ such that $a'(s)<\eps$ for $s>S_{\epsilon}.$ Putting all these facts back in to \eqref{eq24} gives us that $ |\mu_r(s)|\leq Cr^2\eps.$ This finishes the proof of the lemma. \ep

To finish the proof of the the theorem we need another lemma.
\begin{lemma}\label{g}
 Let $g$ be a monotonic increasing function. Let $\mu_r$ and $\theta$ be the functions introduced before. Let as before $\lambda_0\leq\theta(s)\leq\lambda_1$ for all $s.$ Suppose that for every $r\in\mathbb{R}$,
\begin{equation}\label{star}
g(s+r+\mu_r(s))-g(s)-r\sqrt{\theta(s)}\to 0\hspace{0.1in}\mbox{as}\hspace{0.1in} s\to\infty
\end{equation}

\noindent Then for $s\in\mathbb{R},$ the function $g$ can be expressed as
\begin{equation}
\label{eq25} g(s)= \int_{1}^{s}\sqrt{\theta(y)}dy + \widetilde g(s),
\end{equation}
\noindent where $\widetilde g(s)$ satisfies $\tilde g(s+r)-\widetilde g(s)\to 0$ as $s\to\infty,$ for every $r\in\mathbb{R}.$
\end{lemma}

\bp We proceed via a sequence of steps:
 \vskip 0.08 in \noindent {\bf Step I:} We first claim that under the given conditions,
\begin{equation}\label{claim1}
g(s+r+\mu_r(s))-g(s+r)\to 0\quad\mbox{as}\quad s\to\infty,
\end{equation}

\noindent for every $r\in\mathbb{R}$. If this were to be true, then we have
$$ g(s+r)-g(s)-r\sqrt{\theta(s)}= [g(s+r)- g(s+r+\mu_r(s)]+[g(s+r+\mu_r(s))-g(s)- r\sqrt{\theta(s)}].$$
Each of the quantities in the two bracketed terms above tend to $0$ as $s\to\infty$. The first, due to \eqref{claim1}  and the second due to equation \eqref{star}. Therefore we would have as a result
$$g(s+r)-g(s)-r\sqrt{\theta(s)}\to 0\hspace{0.1in}\mbox{as}\hspace{0.1in} s\to\infty.$$
To start with the proof of our first claim, notice that since $g$ is monotone increasing, and $\mu_r(s)\leq 0$ for $r\in\mathbb{R}$ by Lemma \ref{mu}, we have $$\limsup_{s\to\infty}\hspace{0.05in}[g(s+r+\mu_r(s))-g(s+r)]\leq 0.$$
So all we need to do to substantiate the claim is to prove that
$$\liminf_{s\to\infty}\hspace{0.05in} [g(s+r+\mu_r(s))-g(s+r)]\geq 0.$$
This can be seen as follows. Set $\zeta=s+r$. Then $\zeta\to\infty$ as $s\to\infty$. Let $\epsilon>0$ be given. Then, \hspace{0.02in} $\liminf_{s\to\infty}[g(s+r+\mu_r(s))-g(s+r)]=\\$
\begin{eqnarray}
&=& \liminf_{\zeta\to\infty} [g(\zeta+\mu_r(\zeta-r))-g(\zeta)]\nonumber\\
&\ge& \liminf_{\zeta\to\infty} [g(\zeta-Cr^2\eps)-g(\zeta)]     \nonumber\\
&\ge& \liminf_{\zeta\to\infty} [g(\zeta-Cr^2\eps+\mu_{-Cr^2\eps}(\zeta))-g(\zeta)+ Cr^2\eps\sqrt{\theta(\zeta)} - Cr^2\eps\sqrt{\theta(\zeta)}]\nonumber\\
&\ge& \liminf_{\zeta\to\infty} [g(\zeta-Cr^2\eps+\mu_{-Cr^2\eps}(\zeta))-g(\zeta)+Cr^2\eps\sqrt{\theta(\zeta)}]
+ \liminf_{\zeta\to\infty}-Cr^2\eps\sqrt{\theta(\zeta)}\nonumber\\
&\ge& 0 -\sqrt{\lambda_1}Cr^2\epsilon.
\label{eqliminf}
\end{eqnarray}

\noindent An explanation is required here. In the second and third inequalities above, we used the monotonicity of $g$, the the negativity of $\mu_r$ and the bound on $\mu_r$ appearing in \eqref{eq23}. In the last inequality we are applying the hypothesis of the lemma using $-Cr^2\eps$ in place of $r$, along with the fact that $\lambda_0\leq\theta(s)\leq\lambda_1$. Since $\eps>0$ was arbitrary, this proves equation \eqref{claim1}. As was remarked earlier, \eqref{claim1} implies that

\begin{equation}
\label{eq27} g(s+r)-g(s)-r\sqrt{\theta(s)}\to 0\quad\mbox{as}\quad s\to\infty.
\end{equation}

\vspace{0.08 in}

\noindent {\bf Step II:} We make our second claim which is that it follows from equation \eqref{eq27} that for any $\tau\in\mathbb{R}$,
$$\sqrt{\theta(s+\tau)}-\sqrt{\theta(s)}\to 0\hspace{0.1in}\mbox{as}\hspace{0.1in} s\to\infty, $$

\noindent and furthermore, the convergence is uniform on compact sets, with respect to $\tau$. Indeed, let us first replace $s$ with $s+\tau$ and then $r$ by $-\tau$ in equation \eqref{eq27} which is allowed since \eqref{eq27} holds for all $r\in\mathbb{R}.$ As a result, we have that for all $\tau\in\mathbb{R},$
\begin{equation}\label{prelim trans}
g(s)-g(s+\tau)+\tau\sqrt{\theta(s+\tau)}\to 0\quad\mbox{as}\quad s\to\infty.
\end{equation}

\noindent Now applying equation \eqref{eq27} with $\tau$ in place of $r$ and adding to \eqref{prelim trans}, we get that for each $\tau\in\mathbb{R},$
\begin{equation}\label{translate}
\sqrt{\theta(s+\tau)}-\sqrt{\theta(s)}\to 0\hspace{0.1in}\mbox{as}\hspace{0.1in} s\to\infty.
\end{equation}

\noindent We now use the fact that the convergence occurring in equation \eqref{translate} automatically implies uniform convergence on compacts, with respect to $\tau$ (c.f \cite{Seneta}). This proves our second claim. Our third and last step is to show that the function $g$ has the form as stated in the lemma.

\vskip 0.08 in \noindent {\bf Step III:} Define $\widetilde g(s)=g(s)-\int_{1}^{s}\sqrt{\theta(y)}dy.$ Then, for fixed $r\in\mathbb{R},$
\begin{align}
 \widetilde g(s+r)- \widetilde g(s)&= g(s+r)-g(s)-\int_{s}^{s+r}\sqrt{\theta(y)}dy \label{eq:ali1} \\
                &= \hspace{0.05in} g(s+r)-g(s)-r\sqrt{\theta(s+\tau_{r,s})}\hspace{0.1in}\mbox{(for some $\tau_{r,s},$}\hspace{0.1in}|\tau_{r,s}|\leq |r|) \label{eq:ali2} \\
               &= \hspace{0.05in} \left[g(s+r)-g(s)-r\sqrt{\theta(s)}\right]-r\left[\sqrt{\theta(s+\tau_{r,s})}- \sqrt{\theta(s)}\right] \label{eq:ali3}
\end{align}

\noindent We point out that the mean value theorem for integrals was used in \eqref{eq:ali2}. If we now let $s\to\infty$ in \eqref{eq:ali3}, the quantity in the first bracket goes to $0$ by equation \eqref{eq27} and the second bracket tends to $0$ because of \textbf{Step II} proved above (using $|\tau_{r,s}|\leq |r|$ and uniform convergence on compacts). Therefore $$\widetilde g(s+r)-\widetilde g(s)\to 0\hspace{0.1in}\mbox{as}\hspace{0.1in} s\to\infty.$$

\noindent Hence recalling the definition of $\widetilde g$, we have
$$g(s)= \int_{1}^{s}\sqrt{\theta(y)}dy + \widetilde g(s),$$
as required. This finishes the proof of the lemma.\ep

Observe now that the function $f_1$ appearing in equation \eqref{eq20} satisfies the properties stated for $g$ in lemma \ref{g}. Hence applying that lemma we get
\begin{equation}\label{eqf_1}
f_1(s)= \int_{1}^{s}\sqrt{\theta(y)}dy + \kappa(s),
\end{equation}
where $\kappa(s+r)-\kappa(s)\to 0$ for all $r\in\mathbb{R}.$
Such a function $\kappa$ can be written as
\begin{equation}\label{defk_2}
\kappa(s)=\kappa_1(s)+\int_{1}^{s}\kappa_2(y)dy,
\end{equation}
where $\kappa_1$ has a finite limit at infinity and the function $\kappa_2(y)$ tends to 0 as $y$ tends to infinity. For a reference, see \cite{Seneta}. The conditions on $\kappa_1$ and $\kappa_2$ ensure that
\begin{equation}\label{error}
\dfrac{\kappa(s)}{s}\to 0\hspace{0.1in}\mbox{as}\hspace{0.1in} s\to\infty.
\end{equation}

\noindent Next recall that $s=\phi(t)=\int_{1}^{t}\frac{1}{a(y)}dy$, $f_1=f\circ\phi^{-1}$, $\theta=\lambda\circ\phi^{-1}$ and $f(t)=\log M(t)$. Plugging all of this information into \eqref{eqf_1}, we get $$\log M(t)= \int_{1}^{\phi(t)}\sqrt{\lambda\circ\phi^{-1}(y)}dy + \kappa(\phi(t)).$$

\noindent If we change variables, $y=\phi(x)$ in the above integral, we obtain the following.
\begin{equation}\label{prefinal}
\log M(t)= \int_{t_0}^{t}\sqrt{\lambda(x)}\phi'(x)dx + \kappa(\phi(t)),
\end{equation}
\noindent where $t_0 = \phi^{-1}(1).$ Since $\sqrt{\lambda(x)}$ is bounded from below, we have $\int_{t_0}^{t}\sqrt{\lambda(x)}\phi'(x)dx\geq\tilde C\phi(t)$. Therefore for a given $\epsilon>0$ and large $t$,

\begin{equation}\label{finerror}
\dfrac{\kappa(\phi(t))}{\int_{t_0}^{t}\sqrt{\lambda(x)}\phi'(x)dx}\leq\tilde C_1\dfrac{\kappa(\phi(t))}{\phi(t)}< \tilde C_1 \epsilon,
\end{equation}

\noindent the last inequality being true because of equation \eqref{error}. In equation \eqref{prefinal}, we may write the integral from $t_0$ to $t$  as the sum of an integral from $1$ to $t$ plus a bounded part which we may absorb in the error term. Hence we may finally write $\log M(t)$ in the form,

$$ \log M(t)= (1+\rho(t))\int_{1}^{t}\sqrt{\lambda(x)}\phi'(x)dx,\hspace{0.1in}\mbox{where}\hspace{0.1in} \rho(t)\to 0\hspace{0.05in}\mbox{as}\hspace{0.05in} t\to\infty.$$

\noindent Since $\phi'(x)=\dfrac{1}{a(x)}$, this completes the proof of the theorem.

\ep

\noindent \textbf{Remark :}

\noindent As a special case we consider the paraboloid-type domain $$P_{\alpha}= \{(x,Y)\in \mathbb{R}\times\mathbb{R}^{d-1}: \,x > 0, \,|Y|< Ax^{\alpha}\},$$
where $A>0$ and $0<\alpha<1.$ Banuelos and Carroll, in \cite{Ban-Carr} obtained harmonic measure estimates at infinity in $P_{\alpha}.$ As mentioned before, harmonic measure and $M(r)$ are closely related. Roughly speaking, they are reciprocals of one another. We will now illustrate this principle by showing that $P_{\alpha}$ is a cylinder-like domain, apply our theorem to get asymptotic estimates of $M(r),$ and finally compare it with the harmonic measure estimates of Banuelos and Carroll \cite{Ban-Carr}.

\vspace{0.1in}

\noindent First, we show that $P_{\alpha}$ is a cylinder-like domain. Indeed, in the notation of our theorem, choosing $a(t)=At^{\alpha},$ we have,
$$\Gamma_t = \{(x,Y)\in \mathbb{R}\times\mathbb{R}^{d-1}: \, -\frac{t^{1-\alpha}}{2A} < x < \frac{t^{1-\alpha}}{2A},   \, |Y| < (1+\frac{Ax}{t^{1-\alpha}})^{\alpha}\}.$$

\noindent Letting $t\to\infty$, it is easily seen that in the topology of the Hausdorff metric, $\{\Gamma_t\}$ tends to the cylinder $\mathfrak C = \mathbb{R}\times B$ on every compact subset. Here $B$ is the unit ball in $\mathbb{R}^{d-1}.$ Hence, in this case, the family of cylinders $\mathfrak C_t$, consists of the single cylinder $\mathfrak C$. This shows that $P_{\alpha}$ is a cylinder-like domain. Next, notice that $\lambda(t)\equiv\lambda_1,$ the first eigenvalue of the $(d-1)$ dimensional Dirichlet-Laplacian on the unit ball $B.$ Applying our theorem, we obtain,
$$ \log M(t)= (1+\rho(t))\frac{\sqrt\lambda_1}{1-\alpha}t^{1-\alpha} ,$$
where $\rho(t)\to 0$ as $t$ tends to infinity. As remarked in the introduction, a lower bound for harmonic measure is given by the reciprocal of $M(t).$ For the upper bound, we use Proposition $2$ in \cite{Ban-Carr}. Comparing terms, we obtain the same harmonic measure asymptotics as in Theorem $3$ of \cite{Ban-Carr}.

\vspace{0.1in}

\section {Cone-like Domains}

\vspace{0.1in}

\setcounter{equation}{0}

\subsection {Notation and definitions}

\noindent Let $\Omega$ be an unbounded Lipschitz domain in $\mathbb{R}^d.$ In this section we will denote points in $\mathbb{R}^d$ using spherical coordinates: $\xi=(r,\omega)$ where $r= |\xi|$ and $\omega = \frac{\xi}{r}.$ For $t>0$, we set
\begin{equation}\label{0.1a}
\widetilde\Omega_t=\Omega\cap\{\xi: |\xi|\leq t^2 \}
\end{equation}
\begin{equation}\label{0.2a}
\Gamma_t= \frac{\widetilde\Omega_t}{t}.
\end{equation}

A Lipschitz cone will mean a set of the form $(0,\infty)\times D$, where $D$ is a Lipschitz domain properly contained in the unit sphere, i.e., $D\subsetneq\mathbb{S}^{d-1}.$

\begin{definition}\label{charconst}
The characteristic constant $\alpha_1(D)$ of a Lipschitz domain $D\subsetneq\mathbb{S}^{d-1}$ is defined as the positive root of the equation $s^2+(d-2)s-\lambda=0,$ where $\lambda$ is the first eigenvalue of the Laplace-Beltrami operator $\Delta_{\mathbb{S}}$ acting on $D$. If $D$ is on the sphere $|\xi|= t$, let $\widehat D$ be the projection of $D$ on the unit sphere. Then, we define $\alpha_1(\widehat D) = \alpha_1(D).$ If $\Omega$ is a Lipschitz domain and $t>0,$ we define $\alpha_1^{\Omega}(t):= \alpha_1(t)$ to be the characteristic constant of $\Omega\cap\{\xi: |\xi|= t \}.$
\end{definition}

\vspace{0.1in}

\noindent We remark that from the definition \ref{charconst} above, it follows that if $\Omega = (0,\infty)\times D$ is a Lipschitz cone, then $\alpha_1(t)\equiv\alpha_1(D).$ In this section, when we say first eigenvalue, we automatically mean it with respect to the Laplace-Beltrami operator. We next recall the definition of a slowly varying function. A function $B$ defined on the real line is called slowly varying if for every fixed $\mu>0$,
\begin{equation}\label{eq0.3}
\frac{B(t\mu)}{B(t)}\to 1\hspace{0.05in}\mbox{as $t\to\infty.$}
\end{equation}

\noindent Finally, we denote $d_H(A_1,A_2)$ to be the Hausdorff distance between the sets $A_1$ and $A_2.$ Recall the definition of an ULL domain from Definition \ref{Lipschitz}.

\begin{definition}\label{defcones}
We call an ULL domain $\Omega,$ cone-like, if there exists a family of Lipschitz cones $\{\mathfrak G_t\}$ such that for every compact $K,$\\
$$ d_H(\Gamma_t\cap K, \mathfrak G_t\cap K)\to 0\hspace{0.1in}\mbox{as}\hspace{0.1in} t\to\infty.$$

\noindent We also impose the following conditions

\begin{enumerate}
\item The Lipschitz constants of the cones $\mathfrak G_t$ are uniformly bounded above and below.
\item  The characteristic constant $\alpha(t)$ of the cones $\mathfrak G_t$ is a continuous function of $t.$
\item  There exists positive constants $m_0$ and $m_1$ such that for all $t>0$, $m_0\leq\alpha(t)\leq m_1$.
\end{enumerate}
\end{definition}

\vspace{0.1in}

\noindent We note that $\alpha(t)$ appearing in $(2)$ of definition \ref{defcones} has a different meaning than the $\alpha_1$ appearing in definition \ref{charconst}. However, we now show that $\alpha(t)-\alpha_1(t)\to 0$ as $t\to\infty.$ Indeed, let $\{t_n\}$ be any sequence with $t_n\to\infty.$ Then, $\mathfrak G_{t_n}$ is a family of Lipschitz cones whose Lipschitz constants are bounded above and below. Therefore, c.f.\cite{Buc-But}, we can find a subsequence $\{t_{n_k}\}$ and a Lipschitz cone $\mathfrak G$ such that on every compact $K$
$$ d_H(\mathfrak G_{t_{n_k}}\cap K, \mathfrak G\cap K)\to 0,$$

\noindent as $k\to\infty.$ Since $\mathfrak G$ is a Lipschitz cone, we may write $\mathfrak G = (0, \infty)\times D$ for some Lipschitz domain $D\subset \mathbb{S}^{d-1}.$ We note that condition $(3)$ in our definition \ref{defcones} ensures that $D\neq\mathbb{S}^{d-1}.$ Since $\mathfrak G_{t_{n_k}}$ converges to $\mathfrak G$ in the Hausdorff distance, this means that the first eigenvalue of the $\mathfrak G_{t_{n_k}}$ converge to the first eigenvalue of $\mathfrak G$ c.f. \cite{Hen}. This in particular implies that the characteristic constants converge. That is,
\begin{equation}\label{eqdif_0}
\alpha(t_{n_k})\to\tilde\alpha,
\end{equation}

\noindent where $\tilde\alpha$ is the characteristic constant of the cone $\mathfrak G.$ On the other hand by the definition of cone-like domains, we require that $d_H(\Gamma_{t_{n_k}}\cap K,\mathfrak G_{t_{n_k}}\cap K)\to 0$ on every compact $K$. By the triangle inequality this implies also that $d_H(\Gamma_{t_{n_k}}\cap K,\mathfrak G\cap K)\to 0.$ In particular, if we choose $K$ to be the compact set $\overline{\mathfrak G}\cap \{|\xi|=1\}$, then we have $\Gamma_{t_{n_k}}$ converges in the Hausdorff distance to $\mathfrak G\cap \{|\xi|=1\}.$ However, this implies as before that the characteristic constant of $\Gamma_{t_{n_k}}\cap \{|\xi|=1\}$ converge to the characteristic constant of the cone $\mathfrak G.$ But the characteristic constant of $\Gamma_t\cap \{|\xi|=1\}$ is, by definition, $\alpha_1(t)$. Therefore we have
\begin{equation}\label{eqdif_1}
\alpha_1(t_{n_k})\to\tilde\alpha.
\end{equation}

\noindent Combining equations \eqref{eqdif_0} and \eqref{eqdif_1}, we get $\alpha(t_{n_k})-\alpha_1(t_{n_k})\to 0$ as $k\to\infty.$ Since the sequence $\{t_n\}$ was arbitrary we have $\alpha(t)-\alpha_1(t)\to 0$ as $t\to\infty.$ This finishes the proof of the claim. We now proceed to the main theorem of this section.

\subsection{Theorem for cone-like domains}
\setcounter{equation}{5}
\begin{theorem} \label{tMain}

Let $u$ be a positive harmonic function on a cone-like domain $\Omega.$ Suppose that $u=0$ on $\partial\Omega$. Denote $\widetilde M(r)=\sup_{\{\xi\in\Omega : |\xi|=r\}}u(\xi).$ Then, $\widetilde M(r)$ satisfies the following asymptotic formula.

$$ \log\widetilde M(r)=(1+\rho(r)){\int_{e}^{r}\frac{\alpha(y)}{y}dy},$$

\noindent where $\rho$ is a function satisfying $\rho(r)\to 0$ as $r\to\infty.$

\end{theorem}

\bp Fix a positive sequence $\{t_n\}$, $t_n\to\infty$. Let $\Gamma_{t_n}$ and $\widetilde\Omega_{t_n}$ be as defined in equations \eqref{0.1a} and \eqref{0.2a}. By hypothesis we have that on every compact $K$, $d_H(\Gamma_{t_n}\cap K,\mathfrak G_{t_n}\cap K)\to 0\hspace{0.1in}\mbox{as}\hspace{0.1in} n\to\infty.$ Now we use the fact that the set of Lipschitz cones is compact in the topology of Hausdorff distance, c.f. \cite{Buc-But}. Hence there exists a subsequence $\{t_{n_k}\}$ such that $\{\mathfrak G_{t_{n_k}}\}$ converges to an open set $\mathfrak G$ in the Hausdorff distance. It is evident that the limit $\mathfrak G$ is also a Lipschitz cone. So $\mathfrak G$ has the form $(0,\infty)\times D$ for some Lipschitz domain $D\subset\mathbb{S}^{d-1}.$ Furthermore by condition $(3)$ in definition \ref{defcones}, $D\neq\mathbb{S}^{d-1}.$  Applying the triangle inequality we can now conclude that there exists a subsequence $\{t_{n_k}\}$ such that for every compact K
\begin{equation}\label{eqHlimita}
d_H(\Gamma_{t_{n_k}}\cap K, \mathfrak G\cap K)\to 0\hspace{0.1in}\mbox{as $k\to\infty$.}
\end{equation}

\noindent Hence we may pick a unit vector $e$ which is contained in all the $\Gamma_{t_{n_k}}$ as well as in $\mathfrak G.$ To keep the notations simple we will rename the subsequence $\{t_{n_k}\}$ to be $\{t_n\}$ as well as rename the domains $\Gamma_{t_n}$ and $\widetilde\Omega_{t_n}$ to $\Gamma_n$ and $\widetilde\Omega_n.$  We now observe that the mapping $\xi\mapsto t_n\xi$ maps $\overline{\Gamma}_n$ into $\overline{\widetilde\Omega}_n.$ Furthermore, for $\xi\in\Gamma_n$, define $$ v_n(\xi)=\frac{u(t_n\xi)}{\widetilde M(t_n)}.$$ Note that $v_n$ satisfies $v_n(e)\leq1$ and

\begin{equation}\label{harma}\Delta v_n=\dfrac{(t_n)^2}{\widetilde M(t_n)}\Delta u=0,\hspace{0.1in}\mbox{in $\Gamma_n$.}
\end{equation}

\noindent Hence each $v_n$ is a positive harmonic function on $\Gamma_n$ with values $0$ on the lateral part of $\partial\Gamma_n.$ Recall that $v_n(e)\leq 1.$ So, we see that the conditions of Lemma \ref{lemmageneral} are satisfied with $v_n$ in place of $h_n$ and $\Gamma_n$ in place of $\Lambda_n.$ Applying the lemma we deduce that a subsequence, which we still denote $ v_n$, converges to a function $v$ uniformly on compact subsets of the cone ${\mathfrak G}.$ Since $v$ is the uniform limit of harmonic functions, it is harmonic. We now claim that $v=0$ on $\partial {\mathfrak G}.$

\begin{lemma}\label{5a} Let $v$ be a harmonic function on ${\mathfrak G}$ given by
\begin{equation}\label{eqDefva}
v=\lim_{n\to\infty} v_n, \quad v_n(\xi)=\dfrac{u(t_n\xi)}{\widetilde M(t_n)},\quad \mbox{for $\xi\in \Gamma_n$},
\end{equation}
where the convergence is uniform on compact subsets of $\mathfrak G$ and furthermore $v_n$ vanishes on the lateral portion of $\partial\Gamma_n.$ Then $v$ extends to a continuous function on $\overline{{\mathfrak G}}$ with $v|_{\partial {\mathfrak G}}=0$ and for each $r>0,$
\begin{equation}\label{Supb}
\sup_{\{|\xi|=r\}\cap\mathfrak G}{v_n}\to\sup_{\{|\xi|=r\}\cap\mathfrak G}{v}.
\end{equation}
\end{lemma}

The proof of this lemma is similar to the proof of Lemma \ref{lemma5}, with minor modifications. For the sake of completeness, we give the proof here.

\vspace{0.1in}

\bp Fix real numbers $a < b$ and consider the set $A = \{a< |\xi| < b\}\cap\mathfrak G.$ Let $\epsilon > 0$ be given. Now $\overline A$ is a compact set. Therefore by hypothesis there exists $N = N(\epsilon, A)$ such that for $n\geq N$,
\begin{equation}\label{bdrydisb}
d_H(\Gamma_n\cap\overline A, \mathfrak G\cap\overline A)<\epsilon.
\end{equation}

\vspace{0.1in}

\noindent Furthermore, by taking $N$ larger if necessary, we may assume that for $n\geq N$ only the lateral portion $\partial\Gamma_n^L$ of $\partial\Gamma_n$ intersects $\overline A.$ The Lipschitz cone $\mathfrak G$ has the following property : there exists a fixed cone $C(\theta, h),$ of opening $\theta$ and height $h$ such that at every boundary point of $\mathfrak G$ one can draw an exterior cone congruent to $C(\theta, h).$ Since $d_H(\Gamma_n\cap K, \mathfrak C\cap K)\to 0$ on every compact $K$, $\Gamma_n$ also has the same property with the same $C(\theta, h)$ provided $n$ is large. In particular there exists a fixed $r_0 > 0$ independent of $n$ such that for $n\geq N$ and $z\in\partial\Gamma_n^L$, $B_{r_0}(z)\cap\Gamma_n$ can be represented as the region above the graph of some Lipschitz function, with the corresponding Lipschitz constants being uniformly bounded from above and below by positive constants independent of $n.$

\vspace{0.1in}

\noindent Let us now take any $\epsilon,$ $0 < \epsilon < r_0/100,$ with $r_0$ being as in the preceding paragraph. Let $N$ be the corresponding $N(\epsilon, A)$ for which \eqref{bdrydisb} holds. Let $\partial A^L$ be the lateral portion of $\partial A.$ Denote $(\partial A^L)_{\epsilon}= \{x: d(x, \partial A^L) < \epsilon\}.$ Then
\begin{equation}\label{bdrycloseb}
n\geq N\Rightarrow \partial\Gamma_n^L\cap \{a < |\xi| < b\}\subset (\partial A^L)_{\epsilon}.
\end{equation}

\noindent Next cover $(\partial A^L)_{\epsilon}$ by finitely many balls centered at $\partial A^L,$ i.e,
\begin{equation}\label{coverb}
(\partial A^L)_{\epsilon}\subset \bigcup\limits_{i=1}^m{B(z_i, \delta)},
\end{equation}

\vspace{-0.1in}

\noindent where $z_i\in\partial A^L$ and $\delta = r_0/2.$ This is possible by our choice of $\epsilon$ and $\delta.$ Now take any $x\in A\cap B(z_1, \epsilon).$ For $n\geq N$, if $x\notin\Gamma_n,$ extend the definition by setting $v_n(x) = 0.$ Otherwise $x\in\Gamma_n,$ and using $x\in B(z_1, \epsilon)$ along with \eqref{bdrycloseb}, we will get
\begin{equation}\label{dtolateralb}
d(x, \partial\Gamma_n^L) < 2\epsilon.
\end{equation}

\noindent Furthermore $v_n$ is a positive harmonic function and $v_n = 0$ on the Lipschitz portion $ \partial\Gamma_n^L.$ Therefore, by H\"{o}lder estimates near the boundary $\partial\Gamma_n^L$ (c.f \cite{KenigBook}), there exist constants $C>0$ and $0 < \beta\leq 1$ such that
\begin{equation}\label{Holderb}
v_n(x)\leq C\left(\frac{d(x, \partial\Gamma_n^L)}{\delta/2}\right)^{\beta}\sup_{B(z_1,\delta)\cap\Gamma_n}{v_n}.
\end{equation}

\noindent In general the constants $C$ and $\beta$ in equation \eqref{Holderb} will depend on the Lipschitz character of the domains, however by our remark earlier, $\partial\Gamma_n^L$ have Lipschitz constants uniformly bounded from above and below and hence we may choose $C$ and $\beta$ independent of $n.$ Next, we can control the supremum in \eqref{Holderb} by the Carleson estimate c.f \cite{Aikawa},
\begin{equation}\label{Carlesonb}
\sup_{B(z_1,\delta)\cap\Gamma_n}{v_n}\leq C_1v_n(z_{1,\delta}),
\end{equation}
where $ z_{1,\delta}$ is a point in $B(z_1, \delta)\cap\Gamma_n\cap A$ with $|z_1 - z_{1,\delta}|\geq 3\delta/4= 3r_0/8.$ Again, $C_1$ in equation \eqref{Carlesonb} will usually depend on $n$ but using uniform Lipschitz condition of $\Gamma_n$, we may choose it independent of $n.$ Finally let $A_{r_0/8}$ be the compact set $\{x\in A: d(x, \partial A)\geq r_0/8\}.$ Then using the fact that $\{v_n\}$ converges uniformly on compacts to $v,$ we will have for $n\geq N$ that $v_n$ are all uniformly bounded in $A_{r_0/8}$ by a constant $C_2$ that is independent of $n.$ In particular since $z_{1,\delta}\in A_{r_0/8},$ we have for $n\geq N$
\begin{equation}\label{unifbdb}
v_n(z_{1,\delta})\leq C_2.
\end{equation}

\noindent Combining equations \eqref{dtolateralb}, \eqref{Holderb}, \eqref{Carlesonb} and \eqref{unifbdb}, we get
\begin{equation}\label{smallb}
v_n(x)\leq C_3\epsilon^{\beta},
\end{equation}

\noindent for a constant $C_3$ independent of $n.$ A similar estimate as \eqref{smallb} holds for points $\epsilon$ close to $z_i$ in the other balls $B(z_i, \delta)$ for $2\leq i\leq m$ with the same constant $C_3.$ Since $\bigcup\limits_{i=1}^m{B(z_i, \delta)}$ covers $\partial A^L$, we can conclude with the help of \eqref{smallb} that $v = \lim v_n$ satisfies $v = 0$ on $\partial A^L .$ But $\partial A^L = \{a < |\xi| < b\}\cap\partial\mathfrak G$ (by the definition of $A$). Since $a$ and $b$ were arbitrary to start with, we have $v =0$ on $\partial\mathfrak G.$

\vspace{0.1in}

\noindent Moreover it is also true that $v_n$ converges to $v$ uniformly on each section of the form $\{|\xi| = r\}\cap\mathfrak G.$ This is because, for points $x$ on $\{|\xi|= r\}\cap\mathfrak G$ which are close to the boundary $\{|\xi| =r\}\cap\partial\mathfrak G,$ $v_n$ is small by \eqref{smallb} and since we now know $v = 0$ on $\partial\mathfrak G,$ we can apply the same H\"{o}lder estimates as before to show $v(x)$ also satisfies an estimate of the form \eqref{smallb}. But this means we have $v_n$ converging uniformly to $v$ near the boundary $\{|\xi| =r\}\cap\partial\mathfrak G.$ On the other hand, for points well inside $\{|\xi| = r\}\cap\mathfrak G$ ( meaning a compact), we already know that $v_n$ converges uniformly to $v.$ Therefore $v_n$ converges uniformly to $v$ on the whole section $\{|\xi| = r\}\cap\mathfrak G.$ From uniform convergence we get the convergence of the sup norms. That is,

$$\sup_{\{|\xi|=r\}\cap\mathfrak G}{v_n}\to\sup_{\{|\xi|=r\}\cap\mathfrak G}{v}.$$

\noindent This finishes the proof of the lemma.
\ep

Combining equation \eqref{harma} with Lemma~\ref{5a} we deduce that the function $v$ defined by \eqref{eqDefva} satisfies
\begin{equation}\label{eqPropva}
\Delta v=0 \mbox{ in } {\mathfrak G}, \quad v|_{\partial {\mathfrak G}}=0,\quad v>0  \mbox{ in } {\mathfrak G}.
\end{equation}

\noindent Therefore, $v$ can be represented as
\begin{equation}\label{eq16a}
v(\xi)=A_1 r^{\alpha_0}\psi(\omega),\quad\xi\in {\mathfrak G},\quad\xi=(r,\omega),
\end{equation}

\noindent where $\alpha_0$ is the characteristic constant of the Lipschitz cone $\mathfrak G$ and $\psi$ is the eigenfunction of the Laplace-Beltrami operator corresponding to the principal eigenvalue. For the above statement one can refer to \cite{MY}. Let us now go back to the limit in equation \eqref{Supb}. Recalling the definition of $v_n$ we observe that the supremum on the left hand side of equation \eqref{Supb} is $\dfrac{\widetilde M(t_nr)}{\widetilde M(t_n)}.$ From equation \eqref{eq16a} above, we get that the supremum on the right hand side is $Ar^{\alpha_0}$ for a certain positive constant $A.$ Therefore, in the limit we get for every $r>0$,
$$\frac{\widetilde M(t_nr)}{\widetilde M(t_n)}\to Ar^{\alpha_0}\hspace{0.1in}\mbox{as}\hspace{0.1in} n\to\infty .$$

\noindent Plugging in the value 1 for $r$ we infer that $A=1$. We have now obtained the following.
$$\frac{\widetilde M(t_nr)}{\widetilde M(t_n)}\to r^{\alpha_0} \hspace{0.1in}\mbox{as}\hspace{0.1in} n\to\infty .$$

\noindent Recall that $\mathfrak G$ is the limit of $\mathfrak G_{t_n}$ in the Hausdroff metric. So if we denote $\alpha(t_n)$ and $\alpha$ to be the characteristic constants of $\mathfrak G_{t_n}$ and $\mathfrak G$ respectively, then we have that $\alpha(t_n)\to\alpha_0,$ c.f. \cite{Hen}. So we may rewrite the previous limit in the following way:

$$\frac{\widetilde M(t_nr)}{\widetilde M(t_n)r^{\alpha(t_n)}}\to 1 \hspace{0.1in}\mbox{as}\hspace{0.1in} n\to\infty .$$
Since this holds for every subsequence $\{t_n\}$, we can deduce that
\begin{equation}\label{eq17a}
\frac{\widetilde M(tr)}{\widetilde M(t)r^{\alpha(t)}}\to 1 \hspace{0.1in}\mbox{as}\hspace{0.1in} t\to\infty .
\end{equation}

\noindent Taking logarithms we obtain
\begin{equation}\label{eq18a}
\log\widetilde M(tr)-\log\widetilde M(t)-\alpha(t)\log r\to 0 \hspace{0.1in}\mbox{as}\hspace{0.1in} t\to\infty .
\end{equation}

\noindent Introducing the notation $t=e^y$, $r=e^p$, we can write equation \eqref{eq18a} as
\begin{equation}\label{eq19a}
\log\widetilde M(e^{y+p})-\log\widetilde M(e^y)-p\alpha(e^y)\to 0 \hspace{0.1in}\mbox{as}\hspace{0.1in} y\to\infty .
\end{equation}

\noindent Finally setting
\begin{equation}\label{eqdeftheta}
f(y)=\log\widetilde M(e^y), \hspace{0.05in} \theta(y)=\alpha(e^y),\hspace{0.05in}\eqref{eq19a}\hspace{0.05in}\mbox{transforms into}
\end{equation}
\begin{equation}\label{eq20a}
f(y+p)-f(y)-p\theta(y)\to 0 \hspace{0.1in}\mbox{as}\hspace{0.1in} y\to\infty .
\end{equation}

\noindent We will now proceed through a lemma.

\begin{lemma}\label{calclemmacones}
Let $f$ and $\theta$ be the functions introduced in \eqref{eqdeftheta}. Suppose that for every $p>0$, $f(y+p)-f(y)-p\theta(y)\to 0 \hspace{0.1in}\mbox{as}\hspace{0.1in} y\to\infty .$ Then $f$ can be written in the form
\begin{equation}\label{eq21a}
f(y)= {\int_{1}^{y}\theta(s)ds}+g(y), \hspace{0.1in}\mbox{where}\hspace{0.1in} \dfrac{g(y)}{{\int_{1}^{y}\theta(s)ds}}\to 0 \hspace{0.1in}\mbox{as}\hspace{0.1in} y\to\infty.
\end{equation}
\end{lemma}

\bp The proof of this lemma is very similar to Steps {\bf II} and {\bf III} in lemma \ref{g} of Section $2.$ If we use those two Steps, we will get $ f(y)= {\int_{1}^{y}\theta(s)ds}+g(y)$ where $g$ is a function which satisfies
\begin{equation}\label{eq23a}
g(y+p)-g(y)\to 0 \hspace{0.1in}\mbox{as}\hspace{0.1in} y\to\infty.
\end{equation}

\noindent As before, such a function $g$ can be written in the form
\begin{equation}\label{eq24a}
g(y)=\kappa_1(y)+{\int_{1}^{y}\kappa_2(y)dy},
\end{equation}

\noindent where $\kappa_1$ has a finite limit at infinity and $\kappa_2(y)\to 0$ as $y\to\infty$. A proof that $g$ can be written in such a form can be found in \cite{Seneta}. Since $\kappa_2(y)\to 0$ and $\kappa_1$ has a finite limit at infinity, we can read off from \eqref{eq24a} that $\frac{g(y)}{y}\to 0$ as $y\to\infty$. On the other hand $\theta$ is bounded below. Therefore ${\int_{1}^{y}\theta(s)ds}\geq k_0y,$ for some positive constant $k_0$. Combining the previous two statements gives
$$\dfrac{g(y)}{{\int_{1}^{y}\theta(s)ds}}\to 0 $$
\ep

We now recall the definition of $f$ and $\theta$ (from equation \eqref{eqdeftheta}) and plug these into equation \eqref{eq21a}. Then, we get
$$\log\widetilde M(e^y)= {\int_{1}^{y}\alpha(e^s)ds}+g(y).$$

\noindent Changing variables to $z = e^y$, we get $$ \log\widetilde M(z)={\int_{e}^{z}\frac{\alpha(s)}{s}ds}+h(z),$$
where $h(z)= g(\log z)$ and $h$ satisfies
\begin{equation}\label{eqChange}
\dfrac{h(z)}{\int_{e}^{z}\frac{\alpha(s)}{s}ds}\to 0, \hspace{0.1in}\mbox{ as}\hspace{0.1in} z\to\infty.
\end{equation}
Pulling out a factor of ${\int_{e}^{z}\frac{\alpha(s)}{s}ds},$ we get
$$ \log\widetilde M(z) = (1+\rho(z)){\int_{e}^{z}\frac{\alpha(s)}{s}ds},$$ where $\rho(z)= \dfrac{h(z)}{\int_{e}^{z}\frac{\alpha(s)}{s}ds}$ and because of equation \eqref{eqChange}, $\rho(z)\to 0$ as $z\to\infty.$ This concludes the proof of the theorem.

\ep

{\bf Acknowledgements : } The author thanks his thesis advisors Prof. Alexandre Eremenko and Prof. Svitlana Mayboroda for suggesting the problem, for their encouragement, and for the many hours of stimulating discussions.

\bibliographystyle{plain}

\begin{thebibliography}{999}

\bibitem{Aikawa} H.\,Aikawa, {Equivalence between the boundary Harnack principle and the Carleson estimate,} Math. Scand. 103 (2008), no. 1, 61-76.

\bibitem{Ancona} A.\,Ancona, {\it On positive harmonic functions in cones and cylinders,} Rev. Mat. Iberoam. 28 (2012), no. 1, 201–230.

\bibitem{Ban-Carr} R.\,Banuelos, T.\,Carroll, {\it Sharp integrability for Brownian Motion in Parabola-shaped Regions,} J. Funct. Anal. 218 (2005), no. 1, 219-253.

\bibitem{Buc-But} D.\,Bucur, G.\,Buttazo, {\it Variational Methods in Shape Optimization Problems,} Progress in Nonlinear Differential equations and Their Applications, Birkhauser.

\bibitem{Burk} D.\,L.\,Burkholder, {\it Exit times of Brownian Motion, Harmonic Majorization and Hardy Spaces,} Adv. Math. 26 (1977), 182-205.

\bibitem{Car-Hay} T.\,Carroll, W.\,K.\,Hayman, {\it Conformal mapping of parabola-shaped domains,} Comput. Methods Funct. Theory 4 (2004), no. 1, 111-126.

\bibitem{Cran-Li} M.\,Cranston, Y.\,Li, {\it Eigenfunction and Harmonic function estimates in domains with horns and cusps} Comm. Partial Differ. Equ. 22,(1997), 1805-1836

\bibitem{DeBlassie} D.\,DeBlassie, {\it The Martin Kernel for Unbounded Domains,} Potential Anal. 32 (2010), no. 4, 389–404.

\bibitem{Ess-Hal} M.\,Essen, K.\,Haliste, {\it A problem of Burkholder and the existence of harmonic majorants of $|x|^p$ in certain domains in $\mathbb{R}^d.$} Ann. Acad. Sci. Fenn. A.I. Math. 9 (1984) 107-116.

\bibitem{Fr-Hay} S.\,Friedland, W.\,K.\,Hayman, {\it Eigenvalue Inequalities for the Dirichlet Problem on Spheres and the Growth of Subharmonic Functions,} Comment. Math. Helvetici. 51 (1976), 133-61.

\bibitem{GT} D.\,Gilbarg, N.\,Trudinger, {\it Elliptic partial differential equations of second order.} Second edition. Grundlehren der Mathematischen Wissenschaften [Fundamental Principles of Mathematical Sciences], 224. Springer-Verlag, Berlin, (1983).

\bibitem{Hal} K.\,Haliste, {\it Some estimates of harmonic majorants.} Ann. Acad. Sci. Fenn. Ser. A. I. Math., 9 (1984), 117–124.

\bibitem{Hen} A.\,Henrot, {\it Extremum problems for Eigenvalues of Elliptic Operators.} Birkhauser

\bibitem{KenigBook} C.\,Kenig, {\it Harmonic analysis techniques for second order elliptic boundary value problems,} CBMS Regional Conference Series in Mathematics, 83. Published for the Conference Board of the Mathematical Sciences, Washington, DC; by the American Mathematical Society, Providence, RI, 1994.

\bibitem{Miyamoto} I.\,Miyamoto, {\it Harmonic functions in a cylinder which vanish on the boundary,} Japan. J. Math, vol.22, No.2, (1996).

\bibitem{MY} I.\,Miyamoto, H.\,Yoshida, {\it Harmonic functions in a cone which vanish on the boundary,} Math. Nachr. 202 (1999), 177-187.

\bibitem{Seneta} E.\,Seneta, {\it Regularly varying functions, page 5,} Lecture Notes in Mathematics. Published by Springer.

\bibitem{War} S.\,E.\,Warschawaski, {\it On conformal mapping of infinite strips,} Trans. A. M. S 51 (1942) 280-335.


\end{thebibliography}

\end{document}